\documentclass[10pt, reqno]{amsart}

\usepackage[hmargin=0.8in,twosideshift=0in,height=8.6in]{geometry}
\usepackage{amssymb,amsthm}
\usepackage{delarray,verbatim}
\usepackage{natbib}

\usepackage{ifpdf}
\ifpdf
\usepackage[pdftex]{graphicx}
\else
\usepackage[dvips]{graphicx}
\fi

\usepackage{ifpdf}
\usepackage{color}
\definecolor{webgreen}{rgb}{0,.5,0}
\definecolor{webbrown}{rgb}{.8,0,0}
\definecolor{emphcolor}{rgb}{0.95,0.95,0.95}

\usepackage{hyperref}
\hypersetup{%
          colorlinks=true,
          linkcolor=webbrown,
          filecolor=webbrown,
          citecolor=webgreen,
          breaklinks=true}
\ifpdf
\hypersetup{pdftex,
            pdfstartview=FitH, 
            bookmarksopen=true,
            bookmarksnumbered=true
}
\else
\hypersetup{dvips}
\fi

\usepackage{multicol}
\usepackage{multirow}

\numberwithin{equation}{section} 

\newcommand{\F}{\mathcal{F}}

\renewcommand{\P}{\mathbb{P}}
\newcommand{\E}{\mathbb{E}^x}

\newcommand{\eps}{\varepsilon}
\newcommand{\R}{\mathbb{R}}
\renewcommand{\S}{\mathcal{S}}

\renewcommand{\alpha}{r}
\newcommand{\tsigma}{\tilde{\sigma}}

\newtheorem{thm}{Theorem}[section]
\newtheorem{lemma}{Lemma}[section]
\newtheorem{remark}{Remark}[section]
\newtheorem{cor}{Corollary}[section]

\numberwithin{equation}{section}

\title{ On the Perpetual American Put Options for Level Dependent Volatility Models with Jumps}

\author{Erhan Bayraktar }
\address[E. Bayraktar]{Department of
  Mathematics, University of Michigan, Ann Arbor, MI 48109}
\email{erhan@umich.edu}
\thanks{This research  is supported in part by the National Science Foundation. }
\thanks{This article was previously circulated as ``Remarks on the Perpetual American Put Option for Jump Diffusion", see  \cite{bayraktar-2007}.}



\begin{document}

\begin{abstract}
We prove that the perpetual American put option price of level dependent volatility model with compound Poisson jumps is convex and
 is the classical solution of its associated
quasi-variational inequality, that it is $C^2$ except at the
stopping boundary and that it is $C^1$ everywhere (i.e. the smooth
pasting condition always holds). 
\end{abstract}

\maketitle

\section{Introduction}

Let 
$(X^0, B)$, $(\Omega^0,\mathcal{F}^0,\mathbb{P}^0)$, $\{\mathcal{F}_t^0\}$ be the unique weak solution of the stochastic differential equation (see p.300 of \cite{MR1121940})
\begin{equation}\label{eq:diffusion-X0}
dX^0_t=\mu X^0_t dt+\sigma(X_t^0) X_t^0 dB_t, \quad X^0_0=x.
\end{equation}
We will assume that  $x \to \sigma(x)$ is strictly positive. We will also assume that for all $x \in (0,\infty)$ there exists $\eps>0$ such that
\begin{equation}
\int_{x-\eps}^{x+\eps}\frac{1+|\mu| y}{\sigma^2(y) y^2}<\infty.
\end{equation}
Our assumptions on $\sigma$ together with precise description of the process at the boundaries of $(0,\infty)$ (these will be given in the next section) guarantee that (\ref{eq:dyn}) has a unique weak solution thanks to Theorem 5.15 on page 341 of \cite{MR1121940}. We will further assume that $x \to \sigma(x)$ is a continuous function.

Let $(\Omega^1,\mathcal{F}^1,\mathbb{P}^1)$ be a probability space hosting a Poisson random measure
$N$ on $\R_+ \times \R_+$ with mean measure $\lambda \nu(dx) dt$ (in
which $\nu$ is a probability measure on $\R_+$). Let us denote the natural filtration of $\int_{\R_+}z N(dt, dz)$ by$\{\mathcal{F}_t^1\}$. Now consider the product probability space $(\Omega,\F,\P)=(\Omega^0 \times \Omega^1,\F^0 \otimes \F^1  ,\P^0 \times \P^1)$. Let us denote by  $\{\mathcal{F}_t^0\}= \{\mathcal{F}_t \otimes \mathcal{F}_t^1\}$. In this new probability space $(\Omega,\F,\P)$ the Wiener process $B$ and the Poisson random measure $N$ are independent and the process the Markov process defined by $\log(X_t)=\log(X_t^0)+\int_{\R_+}z N(dt, dz)$ is adapted to $\{\F_t\}$. Note that the process $X$ satisfies
\begin{equation}\label{eq:dyn}
dX_t=\mu X_tdt+\sigma(X_t)X_tdB_t+ X_{t-} \int_{\R_+} (z-1) N(dt,dz)
\end{equation}
We will assume that the stock price dynamics is given by $X$.
In this framework, if there is a jump at time $t$, the stock price moves from $X_{t-}$ to $Z  X_t$, in which $Z$'s distribution is given by $\nu$.  $Z$ is a positive random variable and note that when $Z <1$ then the stock price X jumps down when $Z>1$ the stock price jumps up. In the Merton jump diffusion model $Z=\exp(Y)$, in which $Y$ is a Gaussian random variable and $\sigma(x)=\sigma$, for some positive constant $\sigma$.
We will
take $\mu=r+\lambda-\lambda \xi$, in which $\xi=\int_{\R_+}x
v(dx)<\infty$ (a standing assumption) so that $X$ is the price of
a security and the dynamics in (\ref{eq:dyn}) are stated under a
risk neutral measure. Different choices of $\lambda$ and $\xi$ gives different risk neural measures, we assume that these parameters are fixed as a result of a calibration to the historical data.
The value function of the perpetual American put option pricing
problem is
\begin{equation}\label{eq:opt-stop}
V(x):=\sup_{\tau \in \widetilde{\S}}\E\{e^{-\alpha \tau}h(X_{\tau})\},
\end{equation}
in which $h(x)=(K-x)^+$ and $\widetilde{\S}$ is the set of $\{\mathcal{F})_t\}$ stopping times.

We will show that $V$ is convex and that it is the classical solution of the associated
quasi-variational inequality, and that the hitting time of the
interval $(0,l_{\infty})$ is optimal for some $l_{\infty} \in
(0,K)$. Moreover, the value function is in $C^1((0,\infty)) \cap
C^2 ((0,\infty)-\{l_{\infty}\})$ (the smooth pasting condition
holds at $l_{\infty}$). Our result can be seen as an extension of \cite{MR1993268} which showed the convexity and smooth fit properties of the infinite horizon American option problem for a constant elasticity of variance (CEV) model, i.e. $\sigma(x)=\sigma x^{-\gamma}$, $\gamma \in (0,1)$, with no jumps. The value function can not be explicitly obtained as in \cite{MR1993268} because there are jumps in our model and the volatility function $x \to \sigma(x)$ is not specified. We will prove the regularity of the value function $V$ by observing
that it is the limit of a sequence of value functions of
optimal stopping problems for another process that does not jump and coincides with $X$ until the first jump time of $X$. This sequence of functions are defined by iterating a certain functional operator, which maps a certain class of convex functions to a
certain class of smooth functions. 
This sequential approximation technique was used in the context of Bayesian quickest change detection
problems in \cite{bdk05}, \cite{bs} and \cite{sd06}. 
A similar methodology was also employed by \cite{menaldi} which represented the Green functions of the integro-partial differential equations in terms of the Green functions of partial differential equation. The sequential representation of the value function is not only useful for the analysis of the behavior of the value function but also it yields good numerical scheme since the sequence of functions constructed converges to the value function uniformly and exponentially fast.  Other,
somewhat similar, approximation techniques were used to approximate the
optimal stopping problems for \emph{diffusions} (not jump
diffusions), see e.g. \cite{alvarez3} for perpetual optimal
stopping problems with non-smooth pay-off functions, and \cite{carr} for finite time horizon
American put option pricing problems for the geometric Brownian motion.

An alternative to our approach would be to use Theorem 3.1 of \cite{sm} which is a verification theorem for the optimal stopping theorem of Hunt processes. This result can be used to study the smooth pasting principle (see Example 5.3 of \cite{sm}). However this approach relies on being able to determine the
Green function of the underlying process explicitly. On the other hand,
\cite{kyprianu} gave necessary and sufficient conditions for the smooth fit principle
principle is satisfied for the American put option pricing problem
for exponential L\'{e}vy processes generalizing the result of
\cite{mordecki}. However, the results of \cite{kyprianu} can not be applied here in general since unless $\sigma(x) =\sigma$, the process $X$ is not an exponential L\'{e}vy process.
Also, we prove that the value function is the classical solution
of the corresponding quasi-variational inequality and that it is convex, which is not
carried out in \cite{kyprianu}.

The next section prepares the proof of our main result
Theorem~\ref{eq:thm-main}. Here is the outline of our
presentation: First, we will introduce a functional operator $J$,
and define a sequence of convex functions $(v_n(\cdot))_{n\geq 0}$
successively using $J$. Second, we will analyze the properties of
this sequence of functions and its limit $v_{\infty}(\cdot)$. This
turns out to be a fixed point of $J$. Then we will introduce a
family of functional operators $(R_l)_{l \in \R}$, study the
properties of such operators, which can be expressed explicitly using the results from classical diffusion theory. The explicit representation of $R_l$ implies that $R_l
f(\cdot)$ satisfies a quasi-variational inequality for any positive
function $f(\cdot)$. Next, we will show that
 $R_l f(\cdot)=J f (\cdot)$, for
a unique $l=l[f]$, when $f$ is in certain class of convex
functions (which includes $v_n(\cdot)$, $0 \leq n \leq \infty$).
Our main result will follow from observing that
$v_{\infty}(\cdot)=Jv_{\infty}(\cdot)=R_{l[v_{\infty}]}v_{\infty}(\cdot)$ and applying optional sampling theorem.

\section{The Main Result (Theorem~\ref{eq:thm-main}) and its Proof}

We will prepare the proof of our main result, Theorem~\ref{eq:thm-main}, in a sequence of lemmata and corollaries. We need to introduce some notation first.
Let us define an operator $J$ through its action on a test
function $f$ as the value function of the following optimal
stopping problem
\begin{equation}\label{eq:defn-Jf}
J f (x)=\sup_{\tau \in
\S}\E\left\{\int_0^{\tau}e^{-(\alpha+\lambda)t}\lambda \cdot S f
(X_t^0)dt+ e^{-(\alpha+\lambda)\tau} h(X^0_\tau)\right\},
\end{equation}
in which
\begin{equation}\label{eq:defn-Sf}
Sf (x)=\int_{\R_+}f(x z) \nu(dz).
\end{equation}
Here, $X^0=\{X^0_t;t \geq 0\}$ is the solution of \eqref{eq:diffusion-X0}, whose infinitesimal generator is given by
\begin{equation}
\mathcal{A}:=\frac{1}{2}\sigma^2(x)x^2\frac{d^2}{dx^2}+\mu x
\frac{d}{dx},
\end{equation}
and $\S$ is the set of $\{\mathcal{F}^0_t\}$ stopping times.
 Let us denote the increasing and
decreasing fundamental solution of the ordinary second order
differential equation $(\mathcal{A }u) (\cdot)-(\alpha+ \lambda)u(\cdot)=0$ by
$\psi(\cdot)$ and $\varphi(\cdot)$ respectively. Let us denote the
Wronskian of these functions by
\begin{equation}
W(\cdot):=\psi'(\cdot)\varphi(\cdot)-\psi(\cdot)\varphi'(\cdot).
\end{equation}
We will assume that $\infty$ is a natural boundary, which implies that 
\begin{equation}
\lim_{x \rightarrow \infty} \psi(x)=\infty, \;\; \lim_{x \rightarrow \infty}\frac{\psi'(x)}{W'(x)}=\infty, \; \;
\lim_{x \rightarrow \infty} \varphi(x)=0, \;\; \lim_{x \rightarrow \infty}\frac{\varphi'(x)}{W'(x)}=0.
\end{equation}
On the other hand, we will assume that zero is either an exit not entrance boundary (e.g. the CEV model, i.e. when $\sigma(x)=\sigma x^{-\gamma}$, $\gamma \in (0,1)$), which implies that
\begin{equation}\label{eq:assump-varphi1}
\lim_{x \rightarrow 0}\psi(x)=0, \;\; \lim_{x \rightarrow 0}\frac{\psi'(x)}{W'(x)}>0, \;\; \lim_{x \rightarrow 0}\varphi(x)<\infty, \;\;\lim_{x \rightarrow 0}\frac{\varphi'(x)}{W'(x)}=-\infty,
\end{equation} 
or a natural boundary (e.g. the geometric Brownian motion, i.e. when $\sigma(x)=\sigma$)
\begin{equation}\label{eq:assump-varphi2}
\lim_{x \rightarrow 0}\psi(x)=0, \;\; \lim_{x \rightarrow 0}\frac{\psi'(x)}{W'(x)}=0, \;\; \lim_{x \rightarrow 0}\varphi(x)=\infty, \;\;\lim_{x \rightarrow 0}\frac{\varphi'(x)}{W'(x)}=-\infty,
\end{equation} 
see page 19 of \cite{MR1912205}. The next lemma shows that the operator $J$ in (\ref{eq:defn-Jf})
preserves boundedness.
\begin{lemma}
Let $f:\R_+ \rightarrow \R_+$ be a bounded function. Then $J f$ is
also bounded. In fact,
\begin{equation}\label{eq:sup-norm-J-f}
0 \leq \|Jf\|_{\infty}\leq
\|h\|_{\infty}+\frac{\lambda}{\alpha+\lambda}\|f\|_{\infty}.
\end{equation}
\end{lemma}
\begin{proof}
The proof follows directly from (\ref{eq:defn-Jf}).
\end{proof}

Let us define a sequence of functions by
\begin{equation}\label{eq:defn-v-n}
v_0(\cdot)=h(\cdot), \quad v_{n+1}=J v_{n}(\cdot), n \geq 0.
\end{equation}
This sequence of functions is a bounded sequence as the next lemma
shows.
\begin{cor}\label{cor:bounded-seq}
Let $(v_{n})_{n \geq 0}$ be as in (\ref{eq:defn-v-n}). For all $n
\geq 0$,
\begin{equation}\label{eq:sup-norm-v-n}
h(\cdot) \leq v_{n}(\cdot) \leq
\left(1+\frac{\lambda}{\alpha}\right)\|h\|_{\infty}.
\end{equation}
\end{cor}
\begin{proof}
The first inequality follows since it may not be optimal to stop
immediately. Let us prove the second inequality using an induction
argument: Observe that $v_0(\cdot)=h(\cdot)$ satisfies
(\ref{eq:sup-norm-v-n}). Assume (\ref{eq:sup-norm-v-n}) holds for
$n$ and let us show that it holds for when $n$ is replaced by $n+1$. Then using
(\ref{eq:sup-norm-J-f})
\begin{equation}
\|v_{n+1}\|_{\infty}=\|Jv_{n}\|_{\infty}\leq
\|h\|_{\infty}+\frac{\lambda}{\alpha+\lambda}\left(1+\frac{\lambda}{\alpha}\right)\|h\|_{\infty}=\left(1+\frac{\lambda}{\alpha}\right)\|h\|_{\infty}.
\end{equation}
\end{proof}

\begin{lemma}\label{lem:preserve}
The operator $J$ in (\ref{eq:defn-Jf}) preserves order, i.e.
whenever for any $f_1,f_2: \R_+ \rightarrow \R_+$ satisfy
$f_1(\cdot) \leq f_2(\cdot)$, then $Jf_1 (\cdot) \leq
Jf_2(\cdot)$. The operator $J$ also preserves convexity, i.e., if
$f:\R_+ \rightarrow \R_+$ is a convex function, then so is
$Jf(\cdot)$.
\end{lemma}

\begin{proof}
The fact that $J$ preserves order is evident from
(\ref{eq:defn-Jf}).  Let us denote $Y_t=e^{-\mu t} X_t^0$; then $Y$ solves
\begin{equation}
dY_t=Y_t \tsigma(Y_t)dB_t,
\end{equation}
in which $\tsigma(Y_t)=\sigma(e^{\mu t} Y_t)$. 
Let us introduce the operators $K, L$ whose actions on a test function $g$ are given by
\begin{equation}
Kg (y):=\sup_{t \geq 0}Lg(y):=\sup_{t \geq 0}\mathbb{E}^{y}\left[\int_0^{t}e^{-(r+\lambda)u}\lambda Sf(e^{\mu t}Y_u)du+e^{-(r+\lambda)t}g(e^{\mu t}Y_t)\right],
\end{equation}
It follows from arguments similar to those of Theorem 9.4 in \cite{MR1609962} that $Jf(y)=\sup_{n}K^{n}h(y)$. Since the supremum of convex functions is convex it is enough to show that $Lh$ is convex.
This proof will be carried out using the coupling arguments presented in the proof of Theorem 3.1 in \cite{MR1620358}.

Let $0<c<b<a$ and for independent Brownian motions $\alpha$, $\beta$ and $\gamma$ define the processes
\begin{equation}
dA_s=A_s \tsigma(A_s)d \alpha_s, \; A_0=a; \quad dB_s=B_s \tsigma(B_s)d \beta_s, \; B_0=b;\quad dC_s=C_s \tsigma(C_s) d \gamma_s, \; C_0=c,
\end{equation}
which are all martingales since $x \to \sigma(x)$ is bounded.
Let us define $H_a=\inf\{u \geq 0:B_u=A_u\}$, $H_c=\{u\geq 0:B_u=C_u\}$ and $\tau(u)=H_a \wedge H_c \wedge u$. If $\tau(u)=u$, then, 
since $f$ is convex (which implies that $Sf$ is also convex) and $A_u \geq B_u \geq C_u$ 
\begin{equation}\label{eq:conv-SF1}
(A_{u}-C_{u})Sf(e^{\mu u}B_{u}) \leq (B_{u}-C_{u}) Sf(e^{\mu u} A_{u})+(A_{u}-B_{u})Sf(e^{\mu u} C_u).
\end{equation}
If $\tau(u)=H_a$, then $(A_u-C_u)h(B_u)$ has the same law as $(B_u-C_u)h(A_u)$ which implies that
\begin{equation}\label{eq:conv-SF2}
\mathbb{E}\left[(A_u-C_u)Sf(e^{\mu u}B_{u}) 1_{\{\tau=H_a\}}\right]=\mathbb{E}\left[(B_u-C_u) Sf(e^{\mu u} A_{u})1_{\{\tau=H_a\}}\right].
\end{equation}
On the other hand, 
\begin{equation}\label{eq:conv-SF3}
\mathbb{E}\left[(A_u-B_u)Sf(e^{\mu u} C_u)1_{\{\tau=H_a\}}\right]=0.
\end{equation}
Likewise,
\begin{equation}\label{eq:conv-SF4}
\mathbb{E}\left[(A_u-C_u)Sf(e^{\mu u}B_{u}) 1_{\{\tau=H_c\}}\right]=\mathbb{E}\left[(B_u-C_u) Sf(e^{\mu u} A_{u})1_{\{\tau=H_c\}}\right]+\mathbb{E}\left[(A_u-B_u)Sf(e^{\mu u} C_u)1_{\{\tau=H_c\}}\right].
\end{equation}
Thanks to \eqref{eq:conv-SF1}-\eqref{eq:conv-SF4} we have that for all $u \leq t$
\begin{equation}\label{eq:conv-SF5}
\mathbb{E}\left[(A_u-C_u)Sf(e^{\mu u}B_{u}) \right] \leq \mathbb{E}\left[(B_u-C_u) Sf(e^{\mu u} A_{u})\right]+\mathbb{E}\left[(A_u-B_u)Sf(e^{\mu u} C_u)\right].
\end{equation}
Since $A,B,C$ are martingales \eqref{eq:conv-SF5} implies
\begin{equation}\label{eq:conv-SF6}
(a-c)\mathbb{E}\left[Sf(e^{\mu u}B_{u}) \right] \leq (b-c)\mathbb{E}\left[Sf(e^{\mu u} A_{u})\right]+(a-b)\mathbb{E}\left[Sf(e^{\mu u} C_u)\right],
\end{equation}
for all $u \leq t$. Similarly,
\begin{equation}\label{eq:conv-SF7}
(a-c)\mathbb{E}\left[h(e^{\mu t}B_{t}) \right] \leq (b-c)\mathbb{E}\left[h(e^{\mu t} A_{t})\right]+(a-b)\mathbb{E}\left[h(e^{\mu t} C_t)\right].
\end{equation}
Equations \eqref{eq:conv-SF6} and \eqref{eq:conv-SF7} lead to the conclusion that $Lh$ is convex.
\end{proof}

As a corollary of Lemma~\ref{lem:preserve} we can state the
following corollary, whose proof can be carried out by induction.
\begin{cor}\label{cor:increas}
The sequence of functions defined in (\ref{eq:defn-v-n}) is an
increasing sequence of convex functions.
\end{cor}

\begin{remark}\label{rem:v-infty}
Let us define,
\begin{equation}\label{eq:defn-v-infty}
v_{\infty}(\cdot):=\sup_{n \geq 0}v_n(\cdot).
\end{equation}
The function $v_{\infty}(\cdot)$ is well defined as a result of
(\ref{eq:sup-norm-v-n}) and Corollary~\ref{cor:increas}. In fact,
it is positive convex because it is the upper envelope of positive
convex functions and it is bounded by the right-hand-side of
(\ref{eq:sup-norm-v-n}).
\end{remark}

We will study the functions $(v_n(\cdot))_{n\geq 0}$ and $v_{\infty}(\cdot)$
more closely, since their properties will be useful in proving our
main result.

\begin{cor}\label{cor:decreasing}
For each $n$, $v_n(\cdot)$ is a decreasing function on
$[0,\infty)$. The same property holds for $v_{\infty}(\cdot)$.
\end{cor}

\begin{proof}
Any positive convex function on $\R_+$ that is bounded from above is decreasing.
\end{proof}

\begin{remark}\label{rem:at-the-abs-bd}
Since $x=0$ is an absorbing boundary, for any $f:\R_+ \rightarrow
\R_+,
$
\begin{equation}
J f(0)=\sup_{t \in \{0,\infty\}}\int_0^{t}e^{-(\alpha+\lambda)s}
\lambda f(0)ds+e^{-(\lambda+\alpha)t}h(0)=\max\left\{\frac{\lambda
}{\lambda+\alpha}f(0), h(0)\right\}.
\end{equation}

\end{remark}
\begin{remark}\label{rem:sharp-bound}
\,\emph{(Sharper upper bounds and the continuity of the
value function)}. \label{rem:sharp-upper-bound} The upper bound in
Corollary \ref{cor:bounded-seq} can be sharpened using
Corollary~\ref{cor:decreasing} and Remark~\ref{rem:at-the-abs-bd}.
Indeed, we have
\begin{equation}
h(\cdot) \leq v_{n}(\cdot) \leq h(0)=\|h\|_{\infty}=K, \quad
\text{for each $n$}, \quad \text{and} \quad h(\cdot) \leq
v_{\infty}(\cdot) \leq h(0)=\|h\|_{\infty}=K.
\end{equation}
It follows from this observation and
Corollary~\ref{cor:decreasing} that the functions $x \rightarrow v_{n}(x)$, for
every $n$, and $x \rightarrow v_{\infty}(x)$, are continuous at
$x=0$. Since they are convex, these functions are continuous on
$[0,\infty)$.
\end{remark}
\begin{remark}\label{rem:right-der}
The  sequence of functions $(v_n(\cdot))_{n \geq 0}$ and its limit
$v_{\infty}$ satisfy
\begin{equation}
D_+v_n(\cdot) \geq -1 \quad \text{for all $n$ and} \quad
D_{+}v_{\infty}(\cdot) \geq -1,
\end{equation}
in which the function $D_+ f (\cdot)$, is the right derivative of
the function $f(\cdot)$. This follows from the facts that
$v_{n}(0)=v_{\infty}(0)=h(0)=K$, and $v_{\infty}(x) \geq v_{n}(x)
\geq h(x)=(K-x)^+$, for all $x \geq 0$, $n \geq 0$, and that the
functions $v_n(\cdot)$, $n \geq 0$, and $v_{\infty}(\cdot)$, are
convex.
\end{remark}

\begin{lemma}
The function $v_{\infty}(\cdot)$ is the smallest fixed point of
the operator $J$.
\end{lemma}

\begin{proof}
\begin{equation}
\begin{split}
v_{\infty}(x)&=\sup_{n\geq 1}v_{n}(x)=\sup_{n \geq 1}\sup_{\tau
\in \S }\E\left\{\int_0^{\tau}e^{-(\alpha+\lambda)t}\lambda \cdot
S v_n (X_t^0)dt+ e^{-(\alpha+\lambda)\tau} h(X^0_\tau)\right\}
\\&=\sup_{\tau
\in \S }\sup_{n \geq
1}\E\left\{\int_0^{\tau}e^{-(\alpha+\lambda)t}\lambda \cdot S v_n
(X_t^0)dt+ e^{-(\alpha+\lambda)\tau} h(X^0_\tau)\right\}
\\&=\sup_{\tau
\in \S }\E\left\{\int_0^{\tau}e^{-(\alpha+\lambda)t}\lambda \cdot
S (\sup_{n \geq 1} v_n) (X_t^0)dt+ e^{-(\alpha+\lambda)\tau}
h(X^0_\tau)\right\}=J v_{\infty}(x),
\end{split}
\end{equation}
in which last line follows by applying the monotone convergence
theorem twice. If $w:\R_+ \rightarrow \R_+$ is another function
satisfying $w(\cdot)=Jw(\cdot)$, then $w(\cdot)=Jw(\cdot) \geq
h(\cdot)=v_0(\cdot)$. An induction argument yields that $w \geq
v_n(\cdot)$, for all $n \geq 0$, from which the result follows.
\end{proof}

\begin{lemma}\label{lem:fix}
The sequence $\{v_n(\cdot)\}_{n \geq 0}$ converges uniformly to
$v_{\infty}(\cdot)$. In fact, the rate of convergence is
exponential:
\begin{equation}
v_{n}(x) \leq v_{\infty}(x) \leq
v_{n}(x)+\left(\frac{\lambda}{\lambda+\alpha}\right)^n
\|h\|_{\infty}.
\end{equation}
\end{lemma}

\begin{proof}
The first inequality follows from the definition of
$v_{\infty}(\cdot)$. The second inequality can be proved by
induction. The inequality holds when we set $n=0$ by
Remark~\ref{rem:sharp-upper-bound}. Assume that the inequality
holds for $n$. Then
\begin{equation}
\begin{split}
v_{\infty}(x)&=\sup_{\tau \in \S
}\E\left\{\int_0^{\tau}e^{-(\alpha+\lambda)t}\lambda \cdot S
v_{\infty} (X_t^0)dt+ e^{-(\alpha+\lambda)\tau}
h(X^0_\tau)\right\}
\\& \leq \sup_{\tau \in \S
}\E\left\{\int_0^{\tau}e^{-(\alpha+\lambda)t}\lambda \cdot S v_n
(X_t^0)dt+ e^{-(\alpha+\lambda)\tau} h(X^0_\tau)+\int_0^{\infty}dt
\,
e^{-(\lambda+\alpha)t}\lambda\left(\frac{\lambda}{\lambda+\alpha}\right)^n
\|h\|_{\infty}\right\}
\\&=v_{n+1}(x)+\left(\frac{\lambda}{\lambda+\alpha}\right)^{n+1}
\|h\|_{\infty}.
\end{split}
\end{equation}
\end{proof} 
In the next lemma, we will introduce a family of operators whose members map positive functions to solutions of quasi-variational inequalities.
\begin{lemma}\label{lem:r-l}
For any $l \in (0,K)$, let us introduce the operator $R_l$ through
its action on a continuous and bounded test function $f:\R_+ \rightarrow \R_+$
by
\begin{equation}\label{eq:defn-Rl}
R_{l}f(x)=\E\left\{\int_0^{\tau_l}e^{-(\lambda+\alpha)t}\lambda
\cdot S f (X_t^0)dt+ e^{-(\alpha+\lambda)\tau_l}
h(X^0_{\tau_l})\right\},
\end{equation}
in which $\tau_l=\inf\{t \geq 0: X^0_t \leq l\}$. Then
\begin{equation}\label{eq:R-l}
\begin{split}
R_l
f(x)=\left[\psi(x)-\frac{\psi(l)}{\varphi(l)}\varphi(x)\right]
\int_x^{\infty} \frac{2 \lambda \varphi(y)}{y^2\sigma^2(y) W(y)}
Sf(y) dy
+\varphi(x) \int_l^{x} \frac{2 \lambda
\left[\psi(y)-\frac{\psi(l)}{\varphi(l)
}\varphi(y)\right]}{y^2 \sigma^2(y) W(y)} Sf(y) dy+
\varphi(x)\frac{h(l)}{\varphi(l)}, \quad x \geq l,
\end{split}
\end{equation}
and $R_lf(x)=h(x)$, for $x \leq l$.
\end{lemma}

\begin{proof}
Let us define
\begin{equation}\label{eq:R-lr-E}
 R_{l,\rho}f(x)=\E\left\{\int_0^{\tau_l \wedge \tau_\rho}e^{-(\lambda+\alpha)t}\lambda
\cdot S f (X_t^0)dt+ e^{-(\alpha+\lambda)(\tau_l \wedge \tau_\rho) }
h(X^0_{\tau_l \wedge \tau_\rho})\right\},
\end{equation}
in which $\tau_{\rho}:=\inf\{t \geq 0: X^0_t \geq \rho\}$. This expression
satisfies the second-order ordinary differential equation
$\mathcal{A}u(x)-(\alpha+\lambda)u(x)+ \lambda Sf(x)=0$ with
boundary conditions $u(l)=h(l)$ and $u(\rho)=h(\rho)$ and therefore can
be written as
\begin{equation}\label{eq:R-lr}
\begin{split}
R_{l,\rho}
f(x)= \bar{\psi}(x)\int_x^{\rho} \frac{2 \lambda
 \bar{\varphi}(y)}{y^2\sigma^2(y)
W(y)} Sf(y) dy
+ \bar{\varphi}(x) \int_l^{x} \frac{2 \lambda
\bar{\psi}(y)}{y^2 \sigma^2(y) W(y)} Sf(y) dy+
h(l)\frac{\bar{\varphi}(x)}{\bar{\varphi}(l)}+h(\rho) \frac{\bar{\psi}(x)}{\bar{\psi}(\rho)},
\end{split}
\end{equation}
 $x \in [l,\rho]$, in which
 \begin{equation}
 \bar{\varphi}(x)=\varphi(x)-\frac{\varphi(\rho)\psi(x)}{\psi(\rho)}, \quad \bar{\psi}(x)=\psi(x)-\frac{\psi(l)\varphi(x)}{\varphi(l)},
 \end{equation}
see e.g. \cite{karlin-taylor} pages 191-204 and \cite{alvarez2}
page 272. Since $\tau_l \wedge \tau_\rho \uparrow \tau_l$  as $\rho
\rightarrow \infty$ applying monotone and bounded convergence
theorems to (\ref{eq:R-lr-E}) gives $R_{l,\rho}(x) \rightarrow
R_{l}(x)$, as $ \rho \rightarrow \infty$, for all $x \geq 0$. Now
taking the limit of (\ref{eq:R-lr}) we obtain (\ref{eq:R-l}).
\end{proof}

\begin{remark}\label{rem:existance}
For any $l \in (0,K)$, the function $R_lf(\cdot)$ is
differentiable everywhere maybe except at $l$. The left derivative
at $l$, $(R_lf)'(l-)=h'(l)$. On the other hand, the
right-derivative of $R_lf(\cdot)$ at $l$ is
\begin{equation}
(R_lf)'(l+)=\left[\psi'(l)-\frac{\psi(l)}{\varphi(l)}\varphi'(l)\right]
\int_l^{\infty} \frac{2 \lambda \varphi(y)}{y^2 \sigma^2(y) W(y)}
Sf(y) dy +\varphi'(l) \frac{h(l)}{\varphi(l)}.
\end{equation}
The natural question to ask is whether we can find a point $l \in
(0,K) $ such that $R_l'(l+)=R_l'(l-)$, i.e.,
\begin{equation}\label{eq:root}
\left[\psi'(l)\varphi(l)-\psi(l)\varphi'(l)\right] \int_l^{\infty}
\frac{2 \lambda \varphi(y)}{y^2 \sigma^2(y) W(y)} Sf(y)
dy=h'(l)\varphi(l)- \varphi'(l) h(l).
\end{equation}
Since $h(l)=0$ and $h'(l)=0$ for $l>K$ and the left-hand-side is
strictly positive, if a solution exists, it has to be less than
$K$. It follows from Corollary 3.2 in \cite{alvarez2} that
\begin{equation}\label{eq:h-F}
\frac{h'(l)\varphi(l)- \varphi'(l)
h(l)}{\psi'(l)\varphi(l)-\psi(l)\varphi'(l)}=-\int_l^{\infty}\frac{2
\varphi(y)}{y^2 \sigma^2(y)W(y)}F(y)dy,
\end{equation}
in which
\begin{equation}
F(x)=(\mathcal{A}-(\alpha+\lambda))h(x), \quad x \geq 0.
\end{equation}
Therefore (\ref{eq:root}) has a solution if and only if there
exists an $l \in (0,K)$ such that
\begin{equation}\label{eq:finding-optimal-l}
 \int_l^{\infty}
\frac{2 \varphi(y)}{y^2\sigma^2(y) W(y)} (\lambda\cdot
Sf(y)+F(y))dy=0.
\end{equation}
Since $h(x)=h'(x)=0$ for $x>K$, for any $0 \leq f(\cdot) \leq K$
there exists a solution to (\ref{eq:finding-optimal-l}) between
$(0,K)$ if
\begin{equation}\label{eq:int-ineq}
\int_{\eps}^{\infty}\frac{2 \varphi(y)}{ y^2 \sigma(y)^2
W(y)}\left(\lambda K-(\lambda+\alpha) (K-y)1_{\{y < K\} }-\mu y
1_{\{y<K\}}\right)dy<0,
\end{equation}
for some $\eps>0$. Our assumptions in \eqref{eq:assump-varphi1} and \eqref{eq:assump-varphi2} guarantee that \eqref{eq:int-ineq} is satisfied (This can be observed from the formula \eqref{eq:h-F} with the proper choices of $h$ and $F$).
\end{remark}

\begin{lemma}\label{lem:rit-left}
Let $f$ be a convex function and let $D_+f(\cdot)$ be the
right-derivative of $f(\cdot)$. Let $R_lf (\cdot)$ be defined as in Lemma~\ref{lem:r-l}. If $D_+f(\cdot) \geq -1$ and
$\|f\|_{\infty} \leq K$, there exists a unique solution to
\begin{equation}\label{eq:defn-l}
(R_lf)'(l)=h'(l)=-1, \quad l \in (0,K),
\end{equation}
in which $R_l(f)$ is as in (\ref{eq:defn-Rl}).
\begin{equation}\label{eq:text-2}
\text{We will denote the unique solution to (\ref{eq:defn-l}) by
$l[f]$}.
\end{equation}
\end{lemma}

\begin{proof}
Existence of a point $l \in (0,K)$ satisfying (\ref{eq:defn-l})
was pointed out in Remark~\ref{rem:existance}. From the same
Remark and especially (\ref{eq:finding-optimal-l}), the uniqueness
of the solution of (\ref{eq:defn-l}) if we can show the following:
\begin{equation}\label{eq:text}
\text{If for any $x \in (0,K)$\,\, $\lambda \cdot Sf(x)+F(x)=0$,
then $D_+G'(x)<0$},
\end{equation}
in which
\begin{equation}
G(l)= \int_l^{\infty} \frac{2 \varphi(y)}{y^2 \sigma^2(y) W(y)}
(\lambda\cdot Sf(y)+F(y))dy, \quad l \geq 0.
\end{equation}
Indeed if (\ref{eq:text}) is satisfied then $G(\cdot)$ is unimodal
and the maximum of $G(\cdot)$ is attained at either $K$ or at a
point $x \in (0,K)$ satisfying (\ref{eq:text}). One should note
that the right-derivative of $G'$, $D_+G'$ exists since
$\lambda\cdot Sf(y)+F(y)$ is convex.

Now, (\ref{eq:text}) holds if and only if
\begin{equation}\label{eq:ineq}
\lambda D_+ (Sf)(x)+F'(x)>0 \quad \text{or equivalently} \quad
\lambda D_+(Sf)(x)+ r+\lambda-\mu>0, \quad x \in (0,K).
\end{equation}
Since $f$ is bounded and positive convex by assumption, it is decreasing.
Therefore, $D_+f(x) \in [-1,0]$, and this in turn implies that
\begin{equation}\label{eq:Sfp-1}
D_+(Sf)(x)=(S(D_+f))(x) \geq -1.
\end{equation}
The equality can be proved using the dominated convergence theorem,
the inequality is from the assumption that $D_+f(x) \geq -1$. Now,
using (\ref{eq:Sfp-1}), it is easy to observe that (\ref{eq:ineq})
always holds when $\xi>1$, since $\mu=r+\lambda-\lambda \xi$.

We still need to prove the uniqueness when $\xi \leq 1$. This uniqueness holds since in this case we have
\begin{equation}\label{eq:Sfx}
 \lambda Sf(x)+F(x)<0, \quad x \in (0,K),
\end{equation}
and $G(\cdot)$ is unimodal and its maximum is attained at $K$.
Indeed, (\ref{eq:Sfx}) holds if
\begin{equation}
\lambda K -\mu x -(\lambda+\alpha)(K-x)<0, \quad x \in (0,K),
\end{equation}
which is the case since $\mu=r+\lambda-\lambda \xi$ and $\xi<1$.
\end{proof}



\begin{lemma}\label{lemma:quasi-Rf}
Given any convex function satisfying $D_+ f(\cdot) \geq -1$ and
$\|f\|_{\infty} \leq K$ let us define
\begin{equation}\label{eq:defn-Rf}
 (R f)(x):=R_{l[f]}f(x), \quad x\geq 0,
\end{equation}
in which $R_lf(\cdot)$ for any $l \in (0,K)$ is defined in
(\ref{eq:defn-Rl}), and $l[f]$ is defined in
Lemma~\ref{lem:rit-left}. Then the function $Rf(\cdot)$ satisfies
\begin{equation}\label{eq:qsi-1}
(R f)(x)=h(x), \quad
x \in (0,l[f]],
\end{equation}
and
\begin{equation}\label{eq:qsi-2}
(\mathcal{A}-(\alpha+\lambda))R f(x)+ \lambda S f(x) =0, \quad
x \in (l[f],\infty).
\end{equation}
Moreover,
\begin{equation}\label{eq:rit-left}
(Rf)'(l[f]-)=(Rf)'(l[f]+).
\end{equation}
\end{lemma}

\begin{proof}
Equation (\ref{eq:rit-left}) is a consequence of
Lemma~\ref{lem:rit-left}. On the other hand the equalities in
(\ref{eq:qsi-1}) and (\ref{eq:qsi-2}) can be proved using
(\ref{eq:R-l}). 
\end{proof}

\begin{lemma}\label{cor:pre-main}
For every $n$, $0 \leq n \leq \infty$, $v_n(\cdot) \in
C^1(0,\infty) \cap C^2((0,\infty)-\{l_n\})$, in which $(l_n)_{n \in \mathbb{N}}$ is an increasing sequence of functions defined by 
$l_{n+1}:=l[v_n], \; 0\leq n <\infty$. Let
$l_{\infty}:=l[v_{\infty}]$.
(We use (\ref{eq:text-2}) to define these quantities.) 
 Moreover, for each
$0 \leq n < \infty$,
\begin{equation}\label{eq:vn1}
v_{n+1}(x)=h(x), \quad (\mathcal{A}-(\alpha+\lambda)) v_{n+1}(x)+\lambda S
v_n(x) \leq 0, \quad x \in (0,l_{n+1}),
\end{equation}
and
\begin{equation}\label{eq:vn1-2}
v_{n+1}(x)>h(x), \quad (\mathcal{A}-(\alpha+\lambda))
v_{n+1}(x)+\lambda Sv_n(x)=0, \quad x \in (l_{n+1},\infty).
\end{equation}
Furthermore, $v_{\infty}(\cdot)$ satisfies
\begin{equation}\label{eq:qvie-1}
v_{\infty}(x)=h(x), \quad (\mathcal{A}-(\alpha+\lambda))
v_{\infty}(x)+\lambda S v_{\infty}(x) \leq 0, \quad x \in (0,l_{\infty}),
\end{equation}
and
\begin{equation}\label{eq:qvie-2}
v_{\infty}(x)>h(x), \quad (\mathcal{A}-(\alpha+\lambda))
v_{\infty}(x)+ \lambda Sv_{\infty}(x)=0, \quad x \in (l_{\infty},\infty).
\end{equation}
\end{lemma}

\begin{proof}
Recall the definition of $(v_n(\cdot))_{n \in N}$ and
$v_{\infty}(\cdot)$ from (\ref{eq:defn-v-n}) and
(\ref{eq:defn-v-infty}) respectively. From the
Remarks~\ref{rem:sharp-bound} and \ref{rem:right-der} we have
that
\begin{equation}\label{eq:v-nareproper}
\|v_{n}(\cdot)\|_{\infty}\leq K , \quad \text{and} \quad D_{+}v_n(\cdot) \geq -1, \quad 0 \leq n \leq \infty.
\end{equation}
Equation~\eqref{eq:v-nareproper} guarantees that $l_n$ is well defined for all $n$. It follows from \eqref{eq:qsi-1} and the fact that $(v_n)_{n \geq 0}$ is an increasing sequence of functions that $(l_n)_{n \in N}$ is an increasing sequence. Thanks to Lemma~\ref{lemma:quasi-Rf}, $Rv_n$ satisfies \eqref{eq:qsi-1} and \eqref{eq:qsi-2} with $f=v_n$. On the other hand,
 when $x<l_{n+1}$ the inequality in \eqref{eq:vn1} is satisfied thanks to the arguments in the proof of Lemma~\ref{lem:rit-left} (see \eqref{eq:ineq}-\eqref{eq:Sfx} and the accompanying arguments).

 Now as a result of a classical verification theorem, which can be proved by using It\^{o}'s lemma, it follows that $Rv_n=Jv_{n}=v_{n+1}$. This proves \eqref{eq:vn1} and \eqref{eq:vn1-2} except for $v_{n+1}(x)>h(x)$, which follows from the convexity of $v_{n+1}$ and the definition of $l_{n+1}$.

Similarly, $Rv_{\infty}=Jv_{\infty}=v_{\infty}$ and as a result $v_{\infty}$ satisfies \eqref{eq:qvie-1} and \eqref{eq:qvie-2}.
\end{proof}

\begin{thm}\label{eq:thm-main}
Let $V(\cdot)$ be the value function of the perpetual American
option pricing problem in (\ref{eq:opt-stop}) and
$v_{\infty}(\cdot)$ the function defined in
(\ref{eq:defn-v-infty}). Then $V(\cdot)=v_{\infty}(\cdot)$
\begin{equation}
V(x)=\E\left\{e^{-\alpha \tau_{l[\infty]}}
h(X_{\tau_{l[\infty]}})\right\},
\end{equation}
in which $l_{\infty}$ is defined as Lemma~\ref{eq:defn-v-infty}.
The value function, $V(\cdot)$, satisfies the quasi-variational
inequalities (\ref{eq:qvie-1}) and (\ref{eq:qvie-2}) and is convex.
\end{thm}
\begin{proof}
Let us define
\begin{equation}
\tau_x:= \inf\{t \geq 0: X_t \leq l_\infty\},
\end{equation}
and
\begin{equation}
M_t:=e^{-rt}v_{\infty}(X_t).
\end{equation}
Recall that $X$ is the jump diffusion defined in (\ref{eq:dyn}).
It follows from Corollary~\ref{cor:pre-main} and $\|v_{\infty}\|_{\infty} \leq K$ that
$\{M_{t \wedge \tau_x}\}_{t \geq 0}$
is a bounded martingale. 
Using the optional sampling theorem we obtain that
\begin{equation}
v_{\infty}(x)=M_0=\E\left\{M_{\tau_x}\right\}=\E\left\{e^{- r \tau_x}v_{\infty}(X_{\tau_x})\right\}=\E\left\{e^{-r \tau_x} (K-X_{\tau_x})^+\right\} \leq V(x).
\end{equation}
On the other hand, as a result of Lemma~\ref{cor:pre-main} and It\^{o}'s formula for semi-martingales $\{M_{t}\}_{t \geq 0}$ is a positive super-martingale. One should note that although $v_{\infty}$ is not $C^{2}$ everywhere, the It\^{o}'s formula in Theorem 71 of \cite{Protter}  can be applied because
the derivative $v'_{\infty}$ is absolutely continuous.
Applying optional sampling theorem for positive super-martingales we have
\begin{equation}
v_{\infty}(x)=M_0 \geq \E\left\{M_{\tau}\right\}=\E\left\{e^{-r \tau}v_{\infty}(X_{\tau})\right\} \geq \E\left\{e^{-r \tau}(K-X_{\tau})^+\right\},
\end{equation}
therefore $v_{\infty}(x) \geq V(x)$, which implies that $v_{\infty}=V$. As a result $V$ satisfies (\ref{eq:qvie-1}) and (\ref{eq:qvie-2}). The convexity of $V$ follows from Remark~\ref{rem:v-infty}.
\end{proof}

{\small
\bibliographystyle{plain}
\bibliography{references}
}

\end{document}